\documentclass{aptpub}

\usepackage[utf8]{inputenc}
\usepackage[T1]{fontenc}
\usepackage{setspace}
\usepackage{lineno}
\usepackage{subcaption}

\usepackage{pifont}
\usepackage{amssymb, amsmath,amsfonts} 
\usepackage[capitalise]{cleveref}
\usepackage{enumitem}

\usepackage{diagbox}
\usepackage{float}
\usepackage{pgf}
\usepackage{pgfplots}
\usepackage{pgfplotstable}
\usepackage{tikz}
\usetikzlibrary{plotmarks} 
\pgfplotsset{compat=1.6}

\usetikzlibrary{positioning,spy,shapes,shapes.misc,arrows,calc,decorations.markings,external}

\pgfplotsset{graf/.style= {
	width=0.5\columnwidth,
    legend pos = north west,
    no markers, 
    xtick distance = {0.2},
    xmin=0, xmax = 1,
    ylabel= bits,
    xlabel= $\tau$}}

\authornames{ M.~Dias {\it et al}} 
\shorttitle{A method for sampling Bernoulli variables} 

\begin{document}

\title{A method for sampling Bernoulli variables} 

\authorone[Federal University of Campina Grande]{Francisco Marcos de Assis}
\authorone[Federal University of Pernambuco]{Juliana Martins de Assis}
\authorone[Federal University of Campina Grande]{Micael Andrade Dias} 


\addressone{Department of Electical Engineering, Federal University of Campina Grande, Campina Grande -- PB, Postal Code: 58.429-140, Brazil. Email: fmarcos@dee.edu.br} 
\addressone{Department of Statistics, Center of Exact and Natural Sciences, Federal University of Pernambuco, Recife -- PE, Postal Code: 50.740-545, Brazil. Email: juliana@de.ufpe.br} 
\addressone{Department of Electical Engineering, Federal University of Campina Grande, Campina Grande -- PB, Postal Code: 58.429-140, Brazil. Email: micael.souza@ee.ufcg.edu.br}

\begin{abstract}
We introduce new method for generating correlated or uncorrelated Bernoulli random variables by using the binary expansion of a continuous random variable with support on the unit interval. We show that when this variable has a symmetric probability density function around $\frac{1}{2}$, its binary expansion provides equiprobable bits over $\{0,1\}$. In addition we prove that when the random variable is uniformly distributed over $[0,1]$, its binary expansion generates independent Bernoulli random variables. Moreover, we give examples where, by choosing some parameterized nonuniform probability density functions over $[0,1]$, samples of Bernoulli variables with specific correlation values are generated. 
\end{abstract}




\section{Introduction}

The generation of a random sequence of bits has applications in complex numerical simulations
\cite{oliver2013fast, symul2011real}.   
Correlated bit streams are especially necessary in the context of stochastic computing \cite{liu2016synthesis, Chen:2013}, and in cryptography \cite{islam2022using}. 
Distributed simulation of random variables are necessary in quantum computing and theoretical computer science \cite{sudan2019communication}. The simulation of correlated Bernoulli vectors and its posterior summation has also applications in biology, where single nucleotide polymorphisms are involved with the risk of certain diseases \cite{lai2021methods}. 

Both real random processes and deterministic systems may produce random (or pseudorandom) sequences of bits
\cite{moysis2020two}. For example, some papers propose the use of chaotic maps for a pseudorandom bit generation
\cite{han2019image, moysis2020two}. Other papers propose using circuits and semiconductor lasers \cite{oliver2013fast, liu2016synthesis}, or even vacuum field fluctuations of an eletromagnetic field \cite{symul2011real} to generate a random sequence of bits.

Mathematically speaking, a sequence of correlated bits may be undestood as the realization of a multivariate Bernoulli vector. We must notice, however, that given fixed marginal Bernoulli distributions, not all correlation matrices of the corresponding multivariate Bernoulli vector are possible \cite{huber2019}. The general correlation structure between any pair of random variables has been exploited in Fréchet-Hoeffding bounds, where any achievable correlation is a convex combination between these bounds. Specifically, there is a convexity parameter $\lambda\in [0,1]$ which assumes value one when the upper correlation bound is achieved \cite{huber2015}. 

Interestingly, there is an important connection between arbitrary multivariate distributions and multivariate
 Bernoulli distributions. Assume the existence of a certain convexity parameter $\lambda$ between a pair of random
  variables from a random vector with fixed marginals. Then, this value $\lambda$ is possible if and only if 
  there exists the same convexity parameter $\lambda$ for a pair of Bernoulli variables (where the Bernoulli
   marginals have mean 1/2)\cite[Theorem 2, p.604]{huber2015}. Thus, the problem of creating a sequence of
correlated Bernoullivariables is equivalent to the more general problem of creating a random vector from other
 multivariate distributions.


In this paper, we introduce a simple algorithm based on the binary expansion of a continuous random variable with
 support on $[0,1]$ for generating sequences of bits.
 The key idea is rooted on well known arithmetic coding algorithm\cite[p.436]{cover2006} 
 adapted by the authors in \cite{MicaelFMA:2023}.

 First we generate a continuous random variable with support in the interval $[0,1]$, next
  we represent it in binary form (numeral$-2$ base) and take a number of binary digits after the point as the
  Bernoulli sample.  For instance if the trial equals $0.72$ we take this binary form  
   $0.101110\ldots $ and the size $6$ sample is $\{1,0,1,1,1,0\}$.
  
  Despite the fact that our method does not allow for 
 arbitrary correlations between any pair of bits, as will be shown in the simulations, it presents some 
 interesting properties: (i) existence regions of equal (or nearly equal) correlations between different pairs 
 of bits, (ii) low implementation complexity and (iii) parametrization by the distribution. By controlling, 
 for example, the beta distribution parameters $\alpha$ and $\beta$, it is possible to control not only the
  correlations between some pairs of bits, but also if the generated bits are uniformly distributed on $\{0,1\}$ or not.

\section{Preliminaries}

In order to present the method for generating a sequence of correlated random bits, we introduce some definitions. Firstly, let $X$ be a continuous random variable with probability density function $f_X$ and support on $[0,1]$. A binary expansion of $X$ with $n$ bits precision will partition the unit interval in $2^n$ disjoint subsets of same length. If we write the $n$-bit binary expansion of $X\in [0,1]$ as $X = 0.B_1B_2\cdots B_n$ such that $X = \sum_{i=1}^nB_i(\frac{1}{2})^i$, one has that the values of the bits are
\begin{align}
	B_1 &= \begin{cases}
		0,\, \mbox{ if } X<1/2,\\
		1,\, \mbox{ if } X\geq 1/2,
	\end{cases} & 
	B_2 &= \begin{cases}
		0,\, \mbox{ if } X\in[0,1/4)\cup[1/2,3/4),\\
		1,\, \mbox{ if } X\in[1/4,1/2)\cup[3/4, 1],
	\end{cases}
\end{align}
\noindent and for any positive integer $i$ of the binary expansion we have:
\begin{eqnarray}
	B_i = \mathbb{I}\{ X \in \cup _{j=1}^{2^{i-1}}[2^{-i}\cdot (2j-1), 2^{-i}\cdot 2j) \},
\end{eqnarray}
where $\mathbb{I}\{A\}$ is the indicator function of an event $A$, which equals 1 when $A$ occurs and 0 otherwise.

\begin{figure}[!htb]
	\centering
	\includegraphics[width=0.8\textwidth]{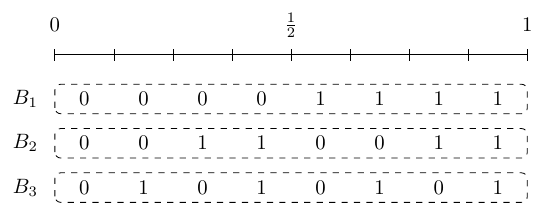}
	\caption{\label{fig:bin-exp-exaple}Unit interval partition according to a $3$-bit binary expansion and the bits corresponding values.}
\end{figure}

In Figure \ref{fig:bin-exp-exaple} we exemplify how the binary expansion works and the corresponding bit values for $n=3$. Each possible $x\in[0,1]$ is contained in one of the partitions and is represented by a unique $n$-bit binary sequence. 

Two questions naturally arise: what is the distribution of $B_i$ and are they correlated? Clearly, $B_i\sim \mbox{Bernoulli}(p_i)$ and $p_i$ and $\mbox{cor}(B_i,B_j)$ depends on the choice of $f_X$, but in which conditions $p_i = \frac{1}{2}$  and $\mbox{cor}(B_i,B_j)=0$? In the next two propositions we prove statements about how $X$ must be distributed in order to give the desired correlated (or uncorrelated) distributions of $B_i$. For the sake of clarity in the following arguments, let us define that a symmetric probability density function of a random variable $X$ around $\frac{1}{2}$ as a  function $f_X:[0,1]\rightarrow \mathbb{R}$ such that for any $\varepsilon\in[0,\frac{1}{2}]$, $f_X(\frac{1}{2}-\varepsilon) = f_X(\frac{1}{2}+\varepsilon)$.
\section{Conditions for Symmetry and Independency}

In this Section we present two propositions about the sequence of bits generated by the binary expansion of the realization of a random variable with support on $[0,1]$. 

\begin{prop}\label{prop:symmetry}
	The bits in the $n$-bit binary expansion of a continuous random variable $X$ with probability density  function $f_X$ and support on the unit interval are $Bern(1/2)$ for any $n\in\mathbb{N}$ if and only if $f_X$ is symmetric around $1/2$.
\end{prop}

\begin{proof}
	($\rightarrow$) Let us define the family $\{A^n\}_{n\geq1}$ of collections of disjoint subsets of $[0,1]$ of length $2^{-n}$:
	\begin{itemize}
		\item $A^1 = \{[0,1/2), [1/2, 1]\}$,
		\item $A^2 = \{[0,1/4),[1/4,1/2),[1/2,3/4), [3/4, 1]\}$
	\end{itemize}
	\noindent and so on. Note that $\bigcup_{A\in A^n}A = [0,1]$ for any $n$. Let us also enumerate the subsets in $A^n$ from $1$ to $n$. Then, the $n$-th bit in the binary expansion of $x$ informs whether $x$ lies in an even of odd numbered subset of $A^n$:
	\begin{itemize}
		\item Odd numbered subset of $A^n$: $b_n = 0$,
		\item Even numbered subset of $A^n$: $b_n = 1$.
	\end{itemize}
	
	Now, define $A_0^n$ and $A_1^n$ as the set of odd and even (respectively) numbered subsets in $A^n$. Then, if $f_X$ is symmetric around $\frac{1}{2}$, 
	\begin{equation}
		\int_{A_0^n[j]}f_X(x)dx = \int_{A_1^n[2^{n-1}-j+1]}f_X(x)dx,\; j=1,\cdots, 2^{n-1}.
	\end{equation}
	
	Also, we have that 
	\begin{align}
		\mbox{Pr}[B_n=0] &= \mbox{Pr}[X\in \bigcup_{j=1}^{2^{n-1}}A_0^n[j]]\\
		&= \sum_{j=1}^{2^{n-1}}\int_{A_0^n[j]}f_X(x)dx\\
		&= \sum_{j=1}^{2^{n-1}}\int_{A_1^n[2^{n-1}-j+1]}f_X(x)dx\\
		&= \mbox{Pr}[X\in A_1^n] = \mbox{Pr}[B_n=1].
	\end{align}
	
	Then, $\mbox{Pr}[B_n=0] = \mbox{Pr}[B_n=1] = \frac{1}{2}$ if $f_X$ is symmetric around $\frac{1}{2}$.
	
	($\leftarrow$) To show that the bits in the binary expansion are $Bern(\frac{1}{2})$ only if $f_X$ is symmetric around $\frac{1}{2}$, assume that $f_X$ is not symmetric around $\frac{1}{2}$. Then, there exists an $\varepsilon$ such that $f_X(\frac{1}{2}-\varepsilon) \neq f_X(\frac{1}{2}+\varepsilon)$. Also, there is some $n\geq 1$ and at least one $j\in\{1,2,\cdots, 2^{n-1}\}$ such that
	\begin{equation}
		\int_{A_0^n[j]}f_X(x)dx \neq \int_{A_1^n[2^{n-1}-j+1]}f_X(x)dx.
	\end{equation}
	
	Then, $\mbox{Pr}[B_n=0] = \mbox{Pr}[B_n=1] = \frac{1}{2}$ only if $f_X$ is symmetric around $\frac{1}{2}$.
\end{proof}

\begin{ex}
	Consider $X$ a random variable that follows the trapezoidal distribution, where $f_X(x)$ is illustrated in Figure \ref{ex:non-sym-dist} and $D-C = \frac{1}{4}$ is fixed. Consider also that $0\leq C \leq \frac{3}{4}$. In the binary expansion of $X$, the probability $Pr[B_1 = 1] = Pr[X\geq \frac{1}{2}] = 1-F_X(\frac{1}{2})$, which is a function of the value of $C$. In \ref{ex:non-sym-prob} we plotted the probabilities $Pr[B_i = 1]$ for the case of a binary expansion with $n=2$. We can see that both bits are equiprobable when $C = \frac{3}{8}$, that is, the $f_X$ is symmetric around $\frac{1}{2}$.
	\begin{figure}[!h]
		\centering
		\begin{subfigure}{0.5\textwidth}
			\centering
			\includegraphics[page=1]{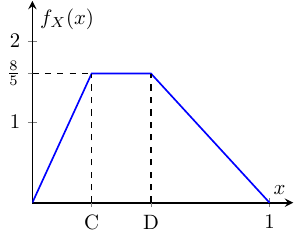}
			\caption{\label{ex:non-sym-dist}}
		\end{subfigure}%
		\begin{subfigure}{0.5\textwidth}
			\centering
			\includegraphics[page=2]{figures/Example01.pdf}
			\caption{\label{ex:non-sym-prob}.}
		\end{subfigure}
		\caption{\label{fig:ex01}(a) Trapezoidal Distribution. (b) Bit probabilities $Pr[B_1 = 1]$ (blue) and $Pr[B_2 = 1]$ (red) in the binary expansion of the trapezoidal distribution as a function of $C$.}
	\end{figure}
\end{ex}

\begin{prop}\label{prop:indep}
	The bits in the $n$-bit binary expansion are independent if $f_X = \mathcal{U}(0,1)$.
\end{prop}

\begin{proof}
	($\rightarrow$) Assume that $f_X = \mathcal{U}(0,1)$. Now consider an $n$-bit binary expansion of $x$ drown from $X$. Each possible $n$-bit sequence corresponds to an interval of length $2^{-n}$. Since $X\sim\mathcal{U}(0,1)$ by hypothesis, the $n$-bit sequences are equiprobable: $p(b_1, \cdots, b_n) = 2^{-n}$. For any two $i,j=1, \cdots, n$, $i\neq j$, we have that
	\begin{equation}
		p(b_i,b_j) = \sum_{\sim b_i,b_j}p(b_1, \cdots,b_n) = (2^{n-2})\cdot2^{-n} = \frac{1}{4}.
	\end{equation}
	
	Then, by factoring $p(b_i,b_j) = p(b_j|b_i)p(b_i)$ and by \ref{prop:symmetry} ensuring that $b_i\sim\operatorname{Bern}(\frac{1}{2})$ for any $i>0$,
	\begin{align}
		\frac{1}{4} = p(b_j|b_i)\frac{1}{2} \rightarrow p(b_j|b_i) = \frac{1}{2} = p(b_j),
	\end{align}
	\noindent from which we conclude that if $X \sim \mathcal{U}(0,1)$ then the bits in the binary expansion are pairwise independent.	
\end{proof}

The other direction of the proposition \ref{prop:indep}, that is, that only an outcome from a uniform distribution could generate a sequence of independent bits in its binary expasion, is not true. We can illustrate this fact with a nonuniform distribution where bits in the binary expansion are also independent. Consider the probability density function given in Equation \ref{eq:exnonuniform}:
\begin{eqnarray}
	f_X(x)=\begin{cases} \dfrac{1}{2}, \mbox{ if }x\in[0,1/2), \\
		\dfrac{3}{2}, \mbox{ if }x\in [1/2, 1] .
	\end{cases} \label{eq:exnonuniform}
\end{eqnarray}


In this case, we may easily find the joint ditribution of the first 3 bits in the binary expansion of $X$, as presented in Tables \ref{tab:b1bj} and \ref{tab:b2b3}, and clearly we can see that they are pairwise independent. 


\begin{table}[H]
\begin{center}
  \begin{tabular}{ |c| c c |c|}
    \hline
    \backslashbox{$B_1$}{$B_j$} & 0  & 1 & $p(b_1)$ \\ \hline
    0 	      & 1/8 & 1/8  &  1/4 \\ 
    1 		  &  3/8  & 3/8 & 3/4   \\ \hline
    $p(b_j)$  & 1/2 & 1/2 & 1 \\ \hline
  \end{tabular}
\caption{Joint distribution of $B_1$ and $B_j$, $j=2,3$, for the binary expansion of an outcome of $X$ with probability density function given in Equation \ref{eq:exnonuniform}.}
\label{tab:b1bj}
\end{center}
\end{table}  

\begin{table}[H]
\begin{center}
  \begin{tabular}{ |c| c c |c|}
    \hline
    \backslashbox{$B_2$}{$B_3$} & 0  & 1 & $p(b_2)$ \\ \hline
    0 	      & 1/4 & 1/4  &  1/2 \\ 
    1 		  & 1/4  & 1/4 & 1/2  \\ \hline
    $p(b_3)$  & 1/2 & 1/2 & 1 \\ \hline
  \end{tabular}
\caption{Joint distribution of $B_2$ and $B_3$, for the binary expansion of an outcome of $X$ with probability density function given in Equation \ref{eq:exnonuniform}.}
\label{tab:b2b3}
\end{center}
\end{table}

As previously shown, when the probability density function $f_X(x)$ is symmetric around $1/2$, than for each individual bit $B_i$ in the binary expansion, we have that $P[B_i=1]=1/2$. The covariance between $B_i$ and $B_j$, $i\neq j$, with $\mathbb{E}[B_i] = \mathbb{E}[B_j] = 1/2$, can be computed as
\begin{align}
	\mbox{cov}(B_i,B_j) &= \mathbb{E}(B_i B_j)-\mathbb{E}(B_i)\mathbb{E}(B_j) \nonumber \\
	&= \mathbb{E}(B_iB_j) - \frac{1}{4}\nonumber \\
	&= Pr[B_i=1, B_j=1] - \frac{1}{4}
\end{align}

It is clear that if $X$ has uniform distribution on $[0,1]$, the generated bits from an outcome of $X$ are independent and thus their correlation is null. In general, for any $f_X(x)$ with support on $[0,1]$ and considering a $3$-bit binary expansion, we have:
\begin{itemize}
	\item $i=1$ and $j=2$:
	\begin{equation}
		Pr[B_i=1, B_j=1] = \int_{\frac{3}{4}}^{1}f_X(x)dx.
	\end{equation}
	\item $i=1$ and $j=3$:
	\begin{equation}
		Pr[B_i=1, B_j=1] = \int_{\frac{5}{8}}^{\frac{6}{8}}f_X(x)dx + \int_{\frac{7}{8}}^{1}f_X(x)dx.
	\end{equation}
	\item $i=2$ and $j=3$:
	\begin{equation}
		Pr[B_i=1, B_j=1] = \int_{\frac{3}{8}}^{\frac{1}{2}}f_X(x)dx + \int_{\frac{7}{8}}^{1}f_X(x)dx.
	\end{equation}
\end{itemize}

When $f_X(x)$ is symmetric around $1/2$, since the bits are $\operatorname{Bern}(\frac{1}{2})$, the correlation between any pair of bits is given by $\mbox{cor}(B_i,B_j) = 4\cdot\mbox{cov}(B_i,B_j)$, which are elements of the correlation matrix $\rho$.

\section{Numerical Results}

\begin{figure}[!tb]
	\centering
	\begin{subfigure}{\textwidth}
		\centering
		\includegraphics[page=1]{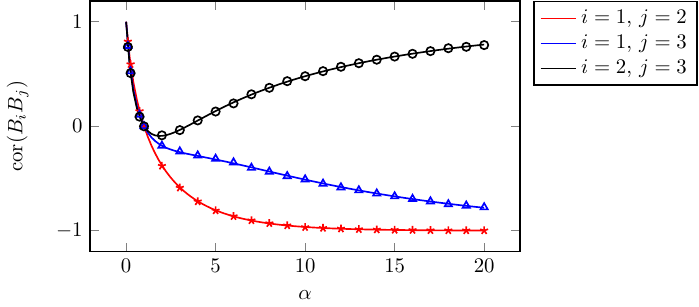}
		\caption{\label{fig:cor-beta}}
	\end{subfigure}
	
	\begin{subfigure}{\textwidth}
		\centering
		\includegraphics[page=2]{figures/Correlations_Plots.pdf}
		\caption{\label{fig:cor-trap}}
	\end{subfigure}
	\caption{\label{fig:corr}Pairwise theoretical (continuous) and estimated (dots, triangles and stars) correlations of bits obtained from binary expansion of (a) beta and (b) trapezoidal random variables. The beta distribution parameters are set to $\alpha=\beta$ and the trapezoidal distribution is symmetric around $\frac{1}{2}$ with upper basis equal to $\Delta$.}
\end{figure}

In Figure \ref{fig:corr} we plotted the theoretical values of correlations obtained by applying the binary expansion to outcomes from beta and trapezoidal distributions with symmetric settings. The symmetry in the distributions is guaranteed by setting $\alpha=\beta$ in the beta distribution, and in the trapezoidal one, $C = \frac{1}{2} - \Delta$ and $D = \frac{1}{2} + \Delta$. We obtained not only the mathematical values of correlation between bits in a 3-bit binary expression, but we also provided estimates of correlation coefficients. Specifically, we simulated $10^5$ outcomes of a beta distribution, for each parameter $\alpha$ in $\{0.1, 0.25, 0.75, 1, 2, 3, \dots, 20\}$. For the trapezoidal distribution we sampled $10^5$ outcomes for each $\Delta$ ranging 10 equally spaced values from 0 to 0.5.

In Figure \ref{fig:corr}, the red curves refer to the correlation between $B_1$ and $B_2$, the blue curves refer to the correlation between $B_1$ and $B_3$ and the black curves refer to the correlation between $B_2$ and $B_3$. We observe from Figure \ref{fig:cor-beta} that setting $\alpha$ with values inside $[0,1]$ results in bits, in a 3-bits binary expasion, that have approximately the same correlation, considering any pair of bits. A similar result occurs when setting $\Delta \in [0.4,0.5]$ in the trapezoidal distribution, as seen in Figure \ref{fig:cor-trap}. 

\section{Concluding remarks}

In this paper we proposed a method for generating random bits, where the bits are a binary expansion of the outcome of a continuous random variable $X$ with support over $[0,1]$. The generated bits may be uniformly distributed over $\{0,1\}$ or not. We proved that the condition of symmetry of the probability density function of $X$ around the value of 1/2 is a necessary and  sufficient condition for uniformly distributed bits. Moreover, we proved that when $X$ is uniformly distributed over $[0,1]$, its outcome binary expansion generate independent bits. Finally, we presented that the method allows the generation of correlated or uncorrelated bits, where the correlation between a pair of bits may be adjusted by parameters of certain distributions. Thus, we believe this method is a promising technique for the generation of a sequence of bits and its many applications.



\fund 
\noindent 
This work was partially supported by Brazilian Council for Scientific and Technological Development (CNPq) under Contract PQ-2 No.311680/2022-4 and the Coordination of Superior Level Staff Improvement (CAPES/PROEX).

\competing 
There were no competing interests to declare which arose during the preparation or publication process of this article.



%
%
%
%

%
%
%
%

\bibliographystyle{APT.bst}
\bibliography{bibliography.bib}

\end{document}